\documentclass[onefignum,onetabnum,a4paper]{article}
\usepackage{epsfig}
\usepackage{epstopdf}
\usepackage{amsfonts}
\usepackage{amssymb}
\usepackage{cases}
\usepackage{gensymb}
\usepackage{graphicx}
\usepackage[T1]{fontenc}

\newcommand{\eref}[1]{(\ref{#1})}
\newcommand{\sref}[1]{section~\ref{#1}}
\newcommand{\fref}[1]{figure~\ref{#1}}

\newcommand{\Fref}[1]{Figure~\ref{#1}}

\title{Chenciner bubbles and torus break-up in a periodically forced delay differential equation} 

\author{{\sc Andrew Keane\footnote{corresponding author:  {\tt a.keane@auckland.ac.nz}} and Bernd Krauskopf} \\ Department of Mathematics, The University of Auckland, \\ Private Bag 92019, Auckland 1142, New Zealand}

\date{February 2018}

\begin{document}

\maketitle

\begin{abstract}

We study a generic model for the interaction of negative delayed feedback and periodic forcing that was first introduced by Ghil \emph{et al.} in the context of the El~Ni\~no Southern Oscillation (ENSO) climate system. This model takes the form of a delay differential equation and has been shown in previous work to be capable of producing complicated dynamics, which is organised by resonances between the external forcing and dynamics induced by feedback.
For certain parameter values, we observe in simulations the sudden disappearance of (two-frequency dynamics on) tori. This can be explained by the folding of invariant tori and their associated resonance tongues. It is known that two smooth tori cannot simply meet and merge; they must actually break up in complicated bifurcation scenarios that are organised within so-called resonance bubbles first studied by Chenciner.

We identify and analyse such a Chenciner bubble in order to understand the dynamics at folds of tori. We conduct a bifurcation analysis of the Chenciner bubble by means of continuation software and dedicated simulations, whereby some bifurcations involve tori and are detected in appropriate two-dimensional projections associated with Poincar\'e sections.
We find close agreement between the observed bifurcation structure in the Chenciner bubble and a previously suggested theoretical picture. As far as we are aware, this is the first time the bifurcation structure associated with a Chenciner bubble has been analysed in a delay differential equation and, in fact, for a flow rather than an explicit map. Following our analysis, we briefly discuss the possible role of folding tori and Chenciner bubbles in the context of tipping.
\end{abstract}

\section{Introduction}
\label{section:intro}

Compared to fold or saddle-node bifurcations of equilibria or periodic solutions, the phenomenon of folding tori (also referred to as a saddle-node bifurcation of tori or a quasiperiodic fold bifurcation) is technically more challenging. In fact, the term ``bifurcation'' is used loosely here, because a fold bifurcation of two smooth invariant tori does not actually exist as such in a generic system. Namely, two smooth invariant tori cannot merge at a single bifurcation point. As the two tori approach the fold locus, they lose normal hyperbolicity (smoothness) and then break up in a bifurcation scenario that involves complicated dynamics \cite{ARO82,BRO98,VIT11}. The bifurcation scenarios that one encounters near the fold are organised in parameter space into so-called resonance bubbles, which were first studied by Chenciner \cite{CHE85a,CHE85b,CHE87b} and are nowadays known as \emph{Chenciner bubbles}.

The loss of smoothness of tori and their subsequent break-up have been observed in various systems; for example, an electric circuit \cite{baptista98}, cardiac fibrillation data \cite{garfinkel97}, the column pendulum \cite{mustafa95} and a particle accelerator model \cite{vrahatis97}. However, understanding the complicated dynamics involved in torus break-up is more challenging. In this paper, we analyse a delay differential equation (DDE), for which a fold bifurcation of tori was observed in \cite{KEA15}. More specifically, we determine precisely how the torus break-up occurs by numerically calculating and identifying the bifurcation structure inside a Chenciner bubble.

The DDE studied in \cite{KEA15} was introduced by Ghil, Zaliapin and Thompson \cite{GHI08} as a conceptual model for the El Ni\~no Southern Oscillation (ENSO) system. ENSO is the climatic driver behind the infamous El Ni\~no events, which are characterised by exceptionally high sea-surface temperatures in the eastern equatorial Pacific Ocean that occur sporadically about every 3--7 years. The DDE model, which we refer to here as the GZT~model, focusses on the interactions of a negative delayed feedback due to ocean-atmosphere coupling with periodic forcing representing the effect of the seasons. It is a special case of a model introduced by Tziperman, Stone, Cane and Jarosh \cite{TZI94} and takes the form:
\begin{equation}
\label{eq:ENSO_model_GHI08}
\dot{{h}}(t) = -{b}\tanh{[{\kappa}{h}(t-\tau_n)]} + {c}\cos{(2\pi t)}.
\end{equation}
The dependent variable $h$ represents the thermocline depth at the eastern boundary of the Pacific Ocean (more specifically, its deviation from the long-term annual mean), which depends on time measured in years.
The parameters $b$ and $c$ of \eref{eq:ENSO_model_GHI08} are amplification factors of the negative delayed feedback and periodic forcing terms, respectively. The delay time $\tau_n$ of the negative feedback mechanism represents the time required for intermediate ocean-atmosphere coupling processes to close the respective feedback loop. In this model the delay is assumed to be constant. The specific form of the coupling function $\tanh(\cdot)$ with coupling strength $\kappa$ is justified in \cite{MUN91}. Throughout, we set $b=1$ and $\kappa=11$ in order to observe the fold bifurcations of tori seen in \cite{KEA15}.
Further details about the climate processes involved can be found in \cite{CLA08}; also see \cite[Sec.~2]{KEA15}.

Since folding tori and Chenciner bubbles are generic or typical phenomena, we expect that the case study presented here will be of interest for the bifurcation analysis of models arising in various applications. 
Because of the simple form that \eref{eq:ENSO_model_GHI08} takes, the GZT model could be viewed more generally as a model for the interaction between negative feedback and periodic forcing; very similar models arise, for example, in human motion control \cite{INS11}, network dynamics \cite{SEM16} and laser systems \cite{SOR15}. Indeed, folding tori could be present whenever the dynamics involves two or more competing frequencies, or put simply, where tori exist; this includes models in the form of ordinary differential equations with a phase space of dimension at least three.

It is the dimensionality of DDEs that allows this seemingly simple model to generate complicated dynamics. Formally, the GZT~model is of the form
\begin{equation}
 \label{eq:DDE_form}
  \dot x(t)=f(t,x(t-\tau_n),\mu)
\end{equation}
where $x \in \mathbb R^n$, $\mu\in\mathbb R^m$ consists of $m$ real parameters and
\begin{equation}
 \label{eq:f_form}
  f: \mathbb R^n \times \mathbb R \times \mathbb R^m \rightarrow \mathbb R^n.
\end{equation}
Note that \eref{eq:DDE_form} is nonautonomous, because $f$ depends on the time $t$ explicitly. For the special case $\tau_n=0$ equation~\eref{eq:DDE_form} is a nonautonomous ordinary differential equation (ODE) and its phase space is $\mathbb R^n \times\mathbb R$, which we also refer to as the physical space. However, as a DDE with $\tau_n>0$, equation~\eref{eq:DDE_form} has the infinite-dimensional space $C([-\tau_n,0];\mathbb R^n) \times \mathbb R$ as its phase space; here $C([-\tau_n,0];\mathbb R^n)$ is the space of continuous functions over the delay interval with values in $\mathbb R^n$. This means that a whole function segment  with values in $\mathbb R^n$ over the time interval $[t_0-\tau_n,t_0]$ is required as part of the initial condition at time $t=t_0$; this is known as an \emph{initial history}. Details on the general theory of DDEs can be found in, for example, \cite{driver77,hale93,STE89}.

In the GZT model \eref{eq:ENSO_model_GHI08} time $t$ enters in the form of a periodic forcing and, hence,  the simplest invariant objects are periodic solutions. The stability of a periodic solution $\Gamma$ is given (as for ODEs) by its Floquet multipliers, which are the eigenvalues of the linearisation along $\Gamma$. For a DDE with a finite number of fixed delays, there exists an infinite number of Floquet multipliers, but there are always at most a finite number of unstable ones \cite{driver77,hale93}. If all the Floquet multipliers (except for the trivial Floquet multiplier 1) lie within the unit circle of the complex plane, then $\Gamma$ is stable; on the other hand, if one or more Floquet multipliers lie outside the unit circle, then $\Gamma$ is of saddle type; together with the trivial Floquet multiplier 1, the corresponding eigendirections span the (finite-dimensional) linear unstable manifold of $\Gamma$. It is convenient for the later discussion to distinguish saddle periodic solutions with one and two unstable Floquet multipliers, and we refer to them as 1-saddle and 2-saddle periodic solutions, respectively. 

Simulations of \eref{eq:ENSO_model_GHI08} can readily be performed with a fixed-step numerical integration method, where the stepsize is given by uniform discretisation of the history segment. Simulations are only able to reveal stable solutions, and saddle solutions must and can be found with continuation software. 
We use the continuation package DDE-Biftool \cite{ENG00,SIE14}, which runs in Matlab and is specifically designed for DDEs. Due to the periodic forcing of the GZT~model, the simplest invariant solutions are periodic orbits. Once a periodic solution to \eref{eq:ENSO_model_GHI08} is found, we utilise DDE-Biftool to numerically continue (or track) the periodic orbit, regardless of its stability, while varying a system parameter. This is achieved by representing the periodic orbit as the solution of an appropriately discretised boundary value problem. DDE-Biftool determines the stability properties of the periodic orbit by calculating its Floquet multipliers, which are used to identify bifurcations. Codimension-one bifurcations --- namely, saddle-node (or fold), period-doubling and torus (or Neimarck-Sacker) bifurcations --- can in turn be continued in a two-parameter plane by placing appropriate constraints on relevant Floquet multipliers. For background information on continuation methods of DDEs see, for example, \cite{ROO07} and for implementation details specific to the GZT~model, see \cite[Sec.~1]{KEA15}.

In our earlier work \cite{KEA15}, we conducted a detailed bifurcation analysis of the GZT~model \eref{eq:ENSO_model_GHI08} with DDE-Biftool. Generally, the dynamics is driven by two independent mechanisms that induce self-sustaining oscillations: the seasonal forcing and the negative delayed feedback. Their interplay can give rise to two-frequency dynamics on an invariant torus, which may be either unlocked or locked. If the frequencies have an irrational ratio, the trajectory on the torus will never close (i.e. it is \emph{unlocked}) and the solution is called quasiperiodic. If the frequencies have a rational ratio, the trajectory is \emph{locked} and a pair of periodic solutions (usually, a stable and a 1-saddle periodic solution) exist on the torus. In parameter space, the regions where the dynamics is locked to some fixed period are known as \emph{resonance tongues} (also called \emph{Arnold tongues}); these are bounded by saddle-node bifurcations of the respective locked periodic solution.

It was observed in \cite{KEA15} that, when certain parameter are changed, the invariant torus loses stability and disappears in a fashion that is analogous to a fold or saddle-node bifurcation of periodic solutions. More specifically, we showed that for certain ranges of $\tau_n$ values, as the parameter $c$ of \eref{eq:ENSO_model_GHI08} is varied, the observed stable torus approaches a co-existing torus of saddle-type in what appears to be a ``fold'' of tori. 

Our approach to investigating folding tori in the GZT~model is to study the locked solutions, organised into resonance tongues, that exist on the tori as they fold. Clearly, as the tori fold, so too do the associated resonance tongues. Inside each resonance tongue, where it folds, exists a Chenciner bubble \cite{BAE07}: within such a bubble an invariant torus undergoes a transition from stable to saddle-type via a loss of normal hyperbolicity and subsequent break up. Therefore, in order to understand the phenomenon of folding tori it is essential to analyse and understand the dynamics in and around such Chenciner bubbles.

Investigating Chenciner bubbles in the context of climate systems could provide important knowledge towards understanding folding tori as a novel type of climate tipping mechanism. This type of tipping should be expected as a feature of systems with more than one frequency due to external forcing and/or internal feedbacks. In fact, oscillations involving multiple frequencies across a wide range are common in global climate systems \cite{GHI91}. Many climate models, as well as observational data, reveal quasiperiodic dynamics (i.e. dynamics on tori); for example, \cite{gamiz07,hartmann08,saha15,sun15,zuehlsdorff08}. Nonetheless, as far as we are aware, ``multi-frequency climate tipping'' has not been considered in the literature yet.

In this paper we consider Chenciner bubbles of folding resonance tongues in the $(\tau_n,c)$-plane of the GZT~model. We determine their structure by bifurcation analysis and demonstrate the resulting complicated dynamics. To this end, we focus on a ${2\!:\!7}$ resonance tongue. It is the largest one in the region of the $(\tau_n,c)$-plane we consider and will therefore form the largest Chenciner bubble. We employ numerical techniques including the continuation software DDE-Biftool to detail the bifurcation structure of the ${2\!:\!7}$ Chenciner bubble. To our knowledge, the investigation in this paper is also the first time that Chenciner bubbles have been analysed in a DDE. Because it has an infinite-dimensional phase space, in order to interpret the dynamics of the GZT~model we consider appropriate projections of the phase space. As part of the bifurcation structure we uncover in the Chenciner bubble, we detect curves of torus bifurcations, neutral saddles, homoclinic bifurcations of periodic solutions, heteroclinic bifurcations of periodic solutions and fold bifurcations of tori in the GZT DDE model. We find that this structure agrees well with the theoretical structure suggested in \cite{BAE07}. One difference we find is that the torus bifurcation inside the Chenciner bubble we study is subcritical, instead of supercritical as in \cite{BAE07}. We give detailed examples of how we calculate two-dimensional projections of the phase space of the DDE model, related to taking a Poincar\'e section, in order to detect the bifurcations that involve tori. In particular, we compute projections of unstable manifolds to find and illustrate locked tori. Finally, given the additional bifurcation structure present within Chenciner bubbles of folding tori, we briefly discuss the potential role of folding tori as a climate tipping mechanism.

The paper is organised as follows. We set the scene in \sref{section:background} by describing the dynamics within a region of the $(\tau_n,c)$-plane where we expect Chenciner bubbles to exist. In \sref{section:theory} we review some of the theory behind Chenciner bubbles, including the bifurcation structure as suggested in \cite{BAE07}. We then present and discuss in \sref{section:chenciner_gzt} a Chenciner bubble of the GZT~model. In sections~\ref{subsection:T_HoC} and~\ref{subsection:HeC_SNT} we present detailed examples of transitions through certain bifurcations in the Chenciner bubble of the GZT DDE model. Finally, in \sref{section:disc} we discuss how our results relate to climate tipping and point to some interesting directions of future work.

\section{Background on the GZT model}
\label{section:background}

\begin{figure}[t]
  \centering
  \vspace*{1mm}
  \includegraphics{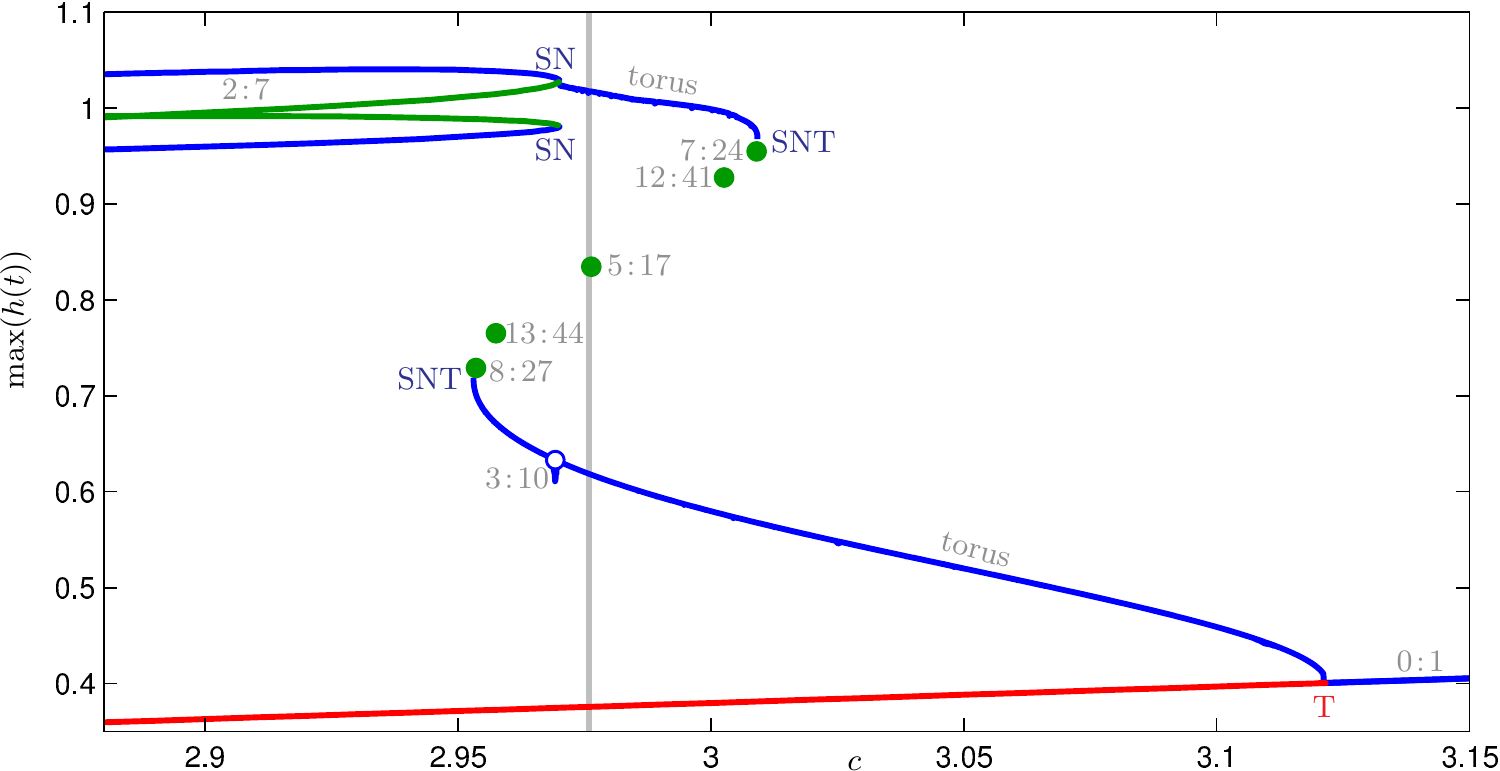}
  \caption{
One-parameter bifurcation diagram of \eref{eq:ENSO_model_GHI08} in $c$ for $\tau_n=0.9395$, where blue, green and red curves indicate stable, 1-saddle and 2-saddle periodic solutions, respectively. Curves of invariant tori are labelled `torus' and resonance tongues are labelled ${p\!:\!q}$. Open blue and filled green circles represent locked solutions that lie on stable and saddle tori, respectively. Torus bifurcations, saddle-node bifurcations of periodic solutions and saddle-node ``bifurcations'' of tori are labelled T, SN and SNT, respectively. The grey vertical line at $c=2.976$ corresponds to the solutions displayed in \fref{fig:SNT_tseries}. Here and throughout $b=1$ and $\kappa=11$.
}
  \label{fig:SNT_1par}
\end{figure}

We now present information on the GZT~model that demonstrates the existence of the phenomenon of folding tori; see also \cite{KEA15}. \Fref{fig:SNT_1par} is a one-parameter bifurcation diagram in $c$ for fixed $\tau_n=0.9395$. The curves labelled `torus' represent maxima of stable tori found by numerical integration, while increasing and decreasing the parameter $c$ slowly enough to ensure that the simulation stays on a particular branch of solutions until stability is lost. As a result, scanning up and down in $c$ reveals bistabilities. Stable, 1-saddle and 2-saddle periodic solutions are represented by their maxima, as computed with the continuation software DDE-Biftool. The point T denotes a torus bifurcations where the ${0\!:\!1}$ periodic solution loses stability while parameter $c$ is decreased. The points SN denote saddle-node bifurcations where stable and 1-saddle periodic solutions of the ${2\!:\!7}$ resonance tongue meet. Because the ${3\!:\!10}$, ${8\!:\!27}$, ${13\!:\!44}$, ${5\!:\!17}$, ${12\!:\!41}$ and ${7\!:\!24}$ resonance tongues are very small, they are represented by circles in \fref{fig:SNT_1par}; open circles represent resonance tongue with both stable and 1-saddle periodic solutions, while filled circles represent resonance tongues with 1-saddle and 2-saddle periodic solutions. In other words, the filled circles represent a branch of solutions on a torus of saddle-type that bridge the two branches of solutions on stable tori. \Fref{fig:SNT_1par} convincingly demonstrates the existence of two loci of folding tori, labelled SNT. These may appear to be just like fold bifurcations of periodic solutions, but as mentioned above, theory predicts more complicated dynamics and bifurcations as the points SNT are approached.

\begin{figure}[t]
  \centering
  \includegraphics[width=\textwidth]{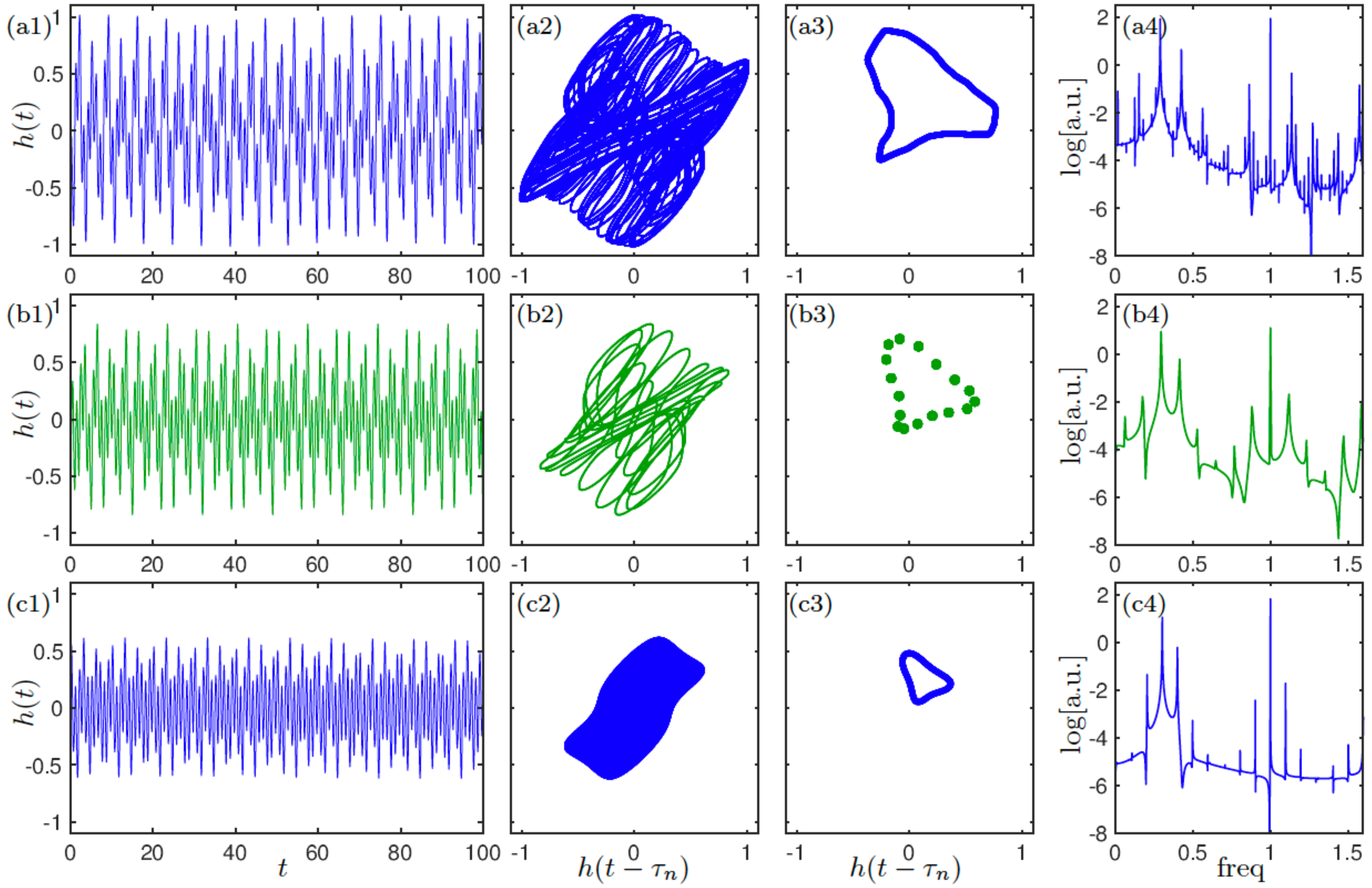}
  \caption{
Solutions of the GZT~model~\eref{eq:ENSO_model_GHI08} for $c=2.976$ and $\tau_n=0.9395$; shown in rows~(a)--(c) are a large stable quasiperiodic solution, a locked saddle periodic solution, and a small stable quasiperiodic solution, respectively. They are displayed as a time series (first column); as a projection in the $(h(t),h(t-\tau_n))$-plane (second column); as a stroboscopic trace in the $(h(t),h(t-\tau_n))$-plane (third column); and as a power spectrum on a logarithmic scale in arbitrary units [a.u.] (fourth column). 
}
  \label{fig:SNT_tseries}
\end{figure}

\Fref{fig:SNT_tseries} shows three co-existing solutions for $c=2.976$ and $\tau_n=0.9395$, as indicated by the vertical line in \fref{fig:SNT_1par}. The solutions in rows~(a) and~(c) of \fref{fig:SNT_tseries} are stable and are found by numerical integration. The locked ${5\!:\!17}$ solution shown in row~(b), on the other hand, is a saddle and calculated with DDE-Biftool. The solutions are displayed (after transients have died down), respectively, as a time series, a projection in the $(h(t),h(t-\tau_n))$-plane, a stroboscopic trace in the $(h(t),h(t-\tau_n))$-plane and a logarithmic power spectrum generated by Fourier transform from the time series over 1000 years. The stroboscopic trace is constructed by plotting after each forcing period the first point of the associated function segment, called the headpoint, in projection onto the $(h(t),h(t-\tau_n))$-plane. Since the forcing period of the GZT~model is 1, this means that in column~3 of \fref{fig:SNT_tseries} we plot $h(t)$ whenever ${t\in\mathbb{N}}$. See \cite{CAL16,KRA03} for a discussion of stroboscopic and Poincar\'e traces of DDEs. 

The solution in row~(a) of \fref{fig:SNT_tseries} belongs to the upper branch of stable tori in \fref{fig:SNT_1par} and it is quasiperiodic (or of very high period). The time series in panel~(a1) of \fref{fig:SNT_tseries} shows a slight modulation of peak heights and the trace of the headpoints in panel~(a3) forms a closed curve. As such, the peaks in the power spectrum in panel~(a4) are incommensurate with the frequency 1 of the seasonal forcing. The solution in row~(b) is a 1-saddle periodic solution from the ${5\!:\!17}$ resonance tongue seen in \fref{fig:SNT_1par}. The phase space projection in panel~(b2) of \fref{fig:SNT_tseries} is a closed curve and the stroboscopic trace in panel~(b3) consists of 17 isolated headpoints. The solution in row~(c) belongs to the lower branch of stable tori in \fref{fig:SNT_1par} and is also quasiperiodic (or of very high period), which is evidenced by the closed loop in the stroboscopic trace in panel~(c3) of \fref{fig:SNT_tseries}. Because the peaks in panel~(a1) are about twice the magnitude of those in panel~(c1), the two stable solutions in rows~(a) and~(c) can be interpreted as two very different climatic states.
 
\begin{figure}[p]
  \centering
  \includegraphics{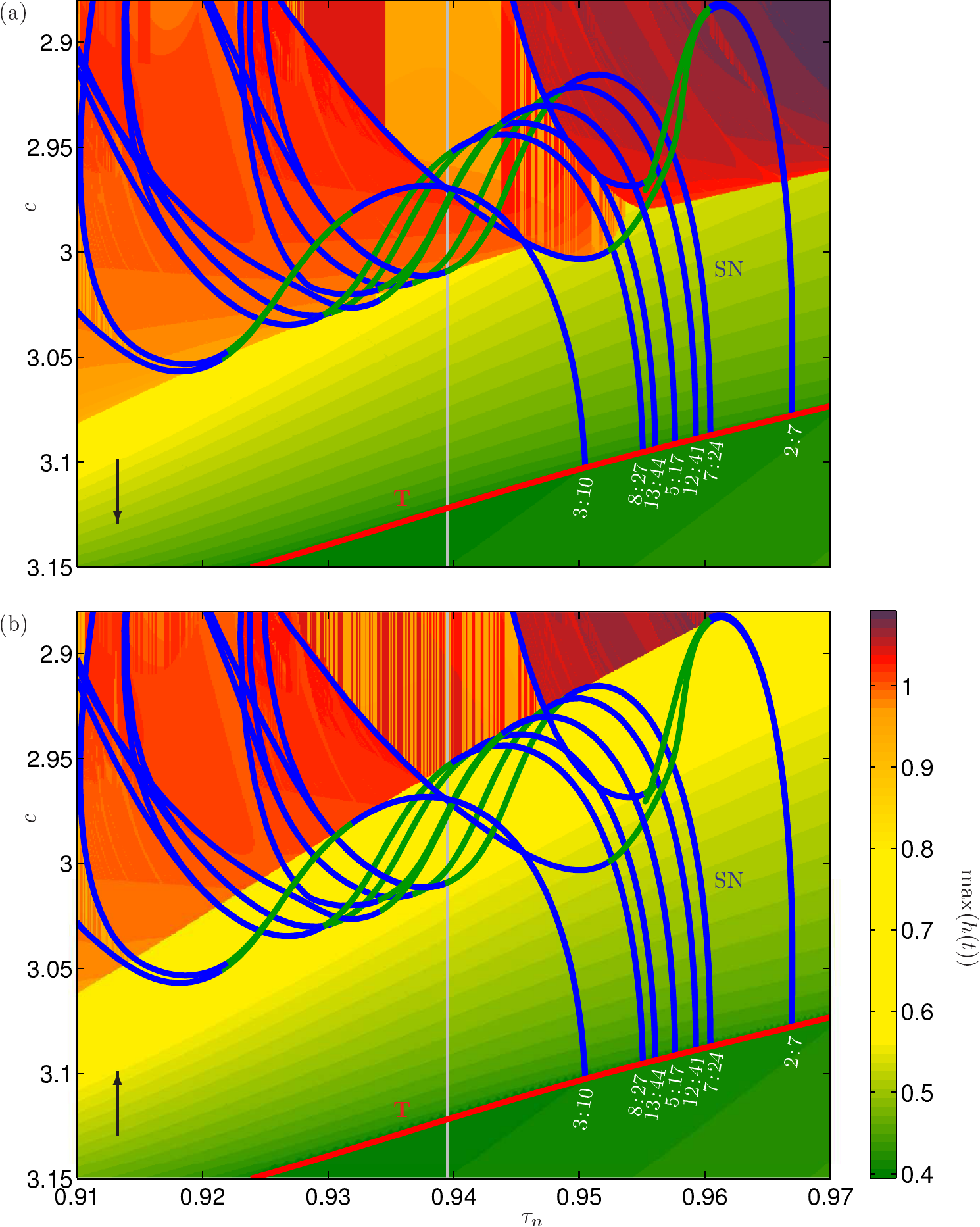}
  \caption{
Maximum maps where maxima of solutions are shown according to the colour bar, for increasing $c$ (a) and for decreasing $c$ (b), as indicated by black arrows. Also shown are the blue/green curves SN of saddle-node bifurcations of $p\!:\!q$ periodic solutions and a red curve T of torus bifurcations. The parts of resonance tongues with green boundaries contain only periodic solutions with at least one unstable Floquet multiplier. The grey vertical line at $\tau_n=0.9395$ indicates the slice shown in \fref{fig:SNT_1par}. 
}
  \label{fig:SNT_example}
\end{figure}

We now illustrate folding tori in the $(\tau_n,c)$-plane. \Fref{fig:SNT_example} shows maximum maps of the GZT~model, which plot the maximum value of attractors found by numerical integration for a range of fixed $\tau_n$ values and for gradually increasing $c$ in panel~(a) and decreasing $c$ in panel~(b), as indicated by the arrows. Maximum maps were calculated in \cite{GHI08,ZAL10}, albeit with a constant fixed initial history for each simulation; in contrast, we scan up and down in parameter $c$ in order to reveal multistabilities, as in \fref{fig:SNT_1par}. We remark that the vertical stripes inside some resonance tongues are not numerical artefacts; they are due to a symmetry property of \eref{eq:ENSO_model_GHI08}, which manifests itself as a bistability in ${p\!:\!q}$ resonance tongues with even $p$ or $q$; see \cite{KEA15} for details. The maximum maps in \fref{fig:SNT_example} are overlaid with curves SN of saddle-node bifurcations of periodic solutions and a curve T of torus bifurcations, computed with DDE-Biftool. The vertical line indicates the parameter slice considered in \fref{fig:SNT_1par}.

The maximum maps in \fref{fig:SNT_example} are divided into two regions of small maxima (green to yellow) and large maxima (orange to red), respectively. The sharp interface in each maximum map, where the maximum of the observed solution changes rapidly from large to small or vice-versa as $c$ is increased or decreased, represents where solutions on tori lose stability and disappear. As demonstrated in \fref{fig:SNT_1par}, this happens where tori fold.

The folding tori are also evident in the bifurcation set of \fref{fig:SNT_example}. Each resonance tongue is bound by two curves SN of saddle-node bifurcations of periodic solutions, emerging from a ${p\!:\!q}$ resonance point along the curve T of torus bifurcations. We observe that the resonance tongues, which are very narrow near the curve T, fold twice in the $c$-direction. Moreover, the foldings of resonance tongues coincide with the sharp interfaces of the maximum maps. Specifically, the two envelopes of folding resonance tongues form the two boundaries where tori fold.

In the parts of the resonance tongue where the boundaries are drawn in blue in \fref{fig:SNT_example}, one finds pairs of stable and 1-saddle periodic solutions; that is, these periodic solutions lie on a stable torus. In the parts of the resonance tongues where the boundaries are drawn in green, on the other hand, the locked periodic solutions have an additional unstable Floquet multiplier (as was checked with DDE-Biftool). Hence, one finds pairs of 1-saddle and 2-saddle periodic solutions; in other words, these locked periodic solutions lie on a torus of saddle-type, with one unstable direction. At $\tau_n=0.9395$ these green segments of resonance tongues correspond to the filled circles in \fref{fig:SNT_1par}. It is where the resonance tongues fold in \fref{fig:SNT_example} and their boundaries change from blue to green that theory, which we briefly review in the next section, predicts the existence of Chenciner bubbles with complicated dynamics.

\section{Theory of bifurcating tori and Chenciner bubbles}
\label{section:theory}

Bifurcations of tori, also referred to as quasiperiodic bifurcations, have been a subject of renewed interest in recent years; for example, see \cite{BAK14,KAM14,KUZ16}. More specifically, quasiperiodic fold bifurcations are studied in, for example, \cite{BRO04,VIT11}, where Chenciner bubbles are identified as regions that could potentially contain interesting dynamics. 

To set the scene for our brief review of relevant theory, consider a family of autonomous vector fields $g$ in three dimensions with two real parameters, $\alpha$ and $\beta$, and a periodic orbit $\Gamma$. For simplicity and to connect with the case considered in this paper, we assume that the vector field is a periodically forced two-dimensional oscillator with a constant forcing period of 1. Hence, the three-dimensional phase space is invariant under translation by 1 of the $t$-axis. Without loss of generality we assume further that the origin is an equilibrium of the oscillator for all $\alpha$ and $\beta$, which generates the periodic orbit $\Gamma$ of period 1 under the forcing. To study the dynamics of the vector field one typically considers the dynamics of the Poincar\'e map $P = P_{\alpha,\beta}$, given as the stroboscopic map of the forcing frequency 1. The origin is a fixed point of $P$ and the (nontrivial) Floquet multipliers of $\Gamma$ are the eigenvalues of the linearisation $DP(0)$ of the map $P$ at the origin. The origin (and hence $\Gamma$) is stable if all the eigenvalues of $DP(0)$ are strictly inside the unit circle. Assume now that the origin is stable for $\alpha <0$ and loses stability at $\alpha=0$, when a pair of complex conjugate eigenvalues of $DP(0)$ moves through the unit circle at $e^{\pm 2\pi i \beta}$. 

While the bifurcation at $\alpha=0$ is of codimension one, it is important to realise that the frequency ratio parameter $\beta$ also plays a role. In fact, it is best to consider the situation in the two-parameter $(\alpha,\beta)$-plane. Either $\beta$ is irrational, that is, $\beta\in\mathbb R\backslash\mathbb Q$, or $\beta$ is rational (or resonant), that is, $\beta=p/q\in\mathbb Q$. In either case an smooth (or normally hyperbolic) invariant circle of $P$ around the origin is born for (sufficiently small) $\alpha > 0$ in what is known as a Neimark-Sacker bifurcation, provided the genericity condition $q\not\in\{1,2,3,4\}$ when $\beta$ is rational.  In the latter case one speaks of a weak $p\!:\!q$ resonance. The cases $q\in\{1,2,3,4\}$ are called strong resonances and they involve other, more complicated bifurcations and unfoldings \cite{arnold77,takens01}. Resonance tongues are found in the $(\alpha,\beta)$-plane, rooted at weak resonance points $\beta=\pm p/q$ along the line $\alpha=0$. For weak resonances, they are bounded by saddle-node bifurcations of a pair of $q$-periodic locked orbits, one of which is an attractor and the other a saddle of $P$. The one-dimensional unstable manifolds of the saddle $q$-periodic orbit meet at the attracting $q$-periodic orbit to form the smooth invariant curve. In between the resonance tongues in the $(\alpha,\beta)$-plane there are smooth curves, rooted at the line $\alpha=0$ at points where $\beta$ is irrational, along which the dynamics on the invariant curve is quasiperiodic, that is, iterates of $P$ fill the invariant curve densely. See, for example, \cite{arnold12,KUZ13,takens01} for further details. 

The above results, including the analysis of the strong resonance, have been obtained by a normal form procedure, developed independently by Arnold \cite{arnold77} and Takens \cite{takens01}, that results in a $\mathbb Z_q$-equivariant planar vector field $G$ whose time-one map is an approximation of the $q$-th iterate of the Poincar\'e  map $P$. This reduces the problem to the analysis of a vector field in two dimensions, instead of three: $q$-periodic points of $P$ are equilibria of $G$, and invariant curves of $P$ are periodic orbits of $G$. 

The road towards a theory of folding tori begins with work by Chenciner in \cite{CHE85a,CHE85b,CHE87b} and, together with Gasull and Llibre, in \cite{CHE87a} on the unfolding of a degenerate Neimark-Sacker bifurcation of codimension two, where the Neimark-Sacker bifurcation changes from being supercritical to subcritical. This means that the bifurcating torus changes from being attracting to being of saddle type. This bifurcation is now often called a Chenciner bifurcation; it necessarily involves a mechanism of merging and annihilation of these two tori and folding resonance tongues. 

The set of further bifurcations associated with a folding $p\!:\!q$ resonance tongue, now referred to as a Chenciner bubble, was considered in \cite{CHE87b,CHE87a}. Derived is a planar vector field of degree two on a cylinder, which corresponds to a fundamental sector of the $\mathbb Z_q$-equivariant planar normal vector field $G$ mentioned above. Effectively, the rotational $\mathbb Z_q$-symmetry has been divided out to obtain a planar vector field on a cylinder, that is, with translational symmetry in the $x$-coordinate on the covering space $\mathbb R^2$. A two-parameter bifurcation diagram is presented  in \cite{CHE87b,CHE87a}. It consists of a pair of saddle-node bifurcations (forming the boundary of the $p\!:\!q$ resonance tongue), which are connected by curves of homoclinic bifurcations and of saddle-node bifurcations of periodic orbits. These additional curves meet at a special point on a line of symmetry where the vector field is Hamiltonian. This bifurcation diagram therefore represents a minimal set of dynamics but not an overall unfolding.

\subsection{Suggested minimal bifurcation diagram of a Chenciner bubble}
\label{subsection:bae07}

In \cite{BAE07} Baesens and MacKay assume an arbitrary weak ${p\!:\!q}$ locking on the torus and that the associated resonance tongue folds in an appropriate parameter plane. They study a more general planar normal form vector field on the cylinder, which includes higher-order terms compared to the one studied by Chenciner. The authors suggest an associated minimal bifurcation diagram for parameter values near the Chenciner bubble, derived with topological argumentations, where the line of Hamiltonian dynamics is unfolded. 

\begin{figure}[t]
  \centering
  \includegraphics{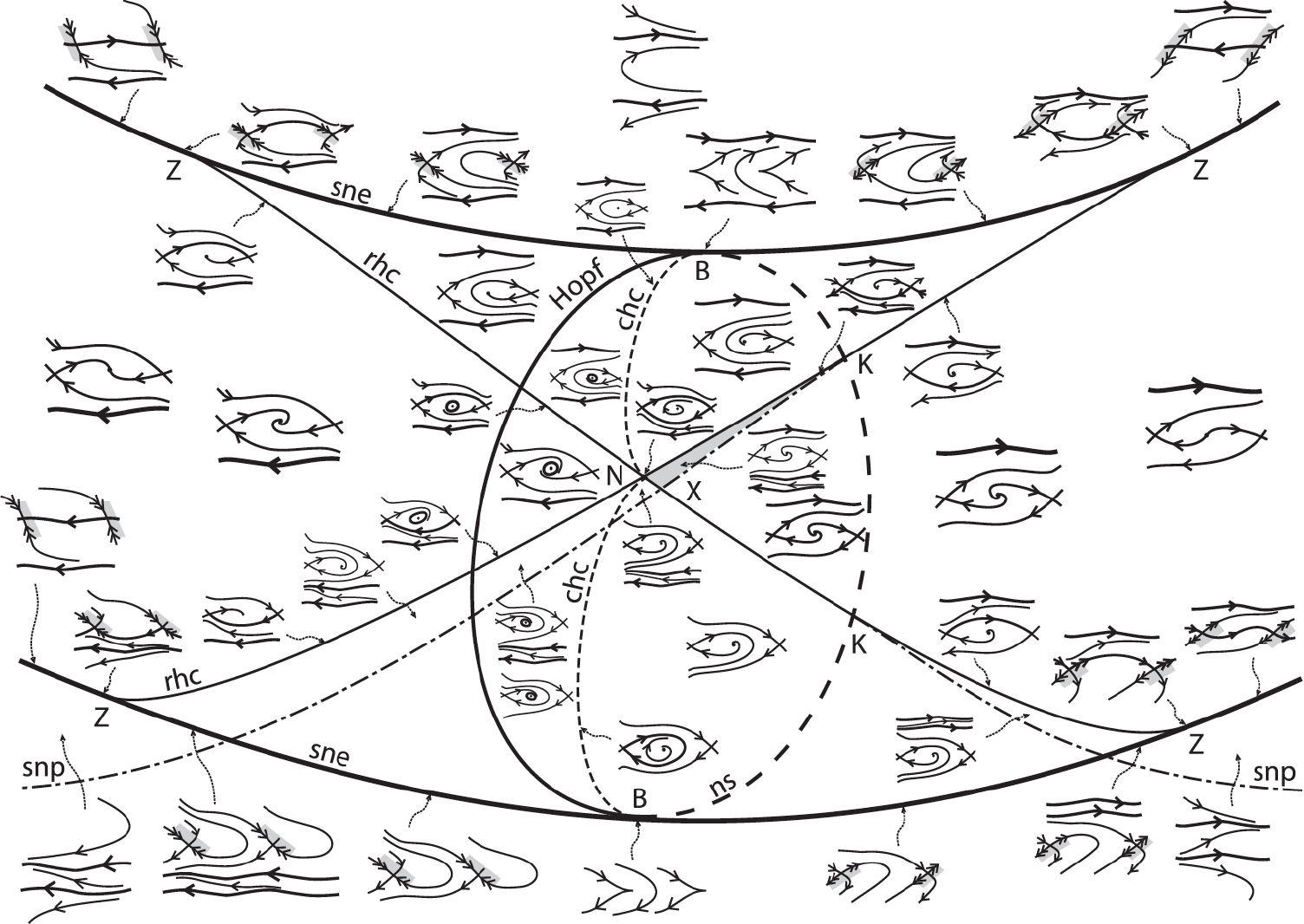}
  \caption{
Bifurcation structure inside a Chenciner bubble with representative phase portraits, as suggested in \cite{BAE07}. Shown are a curve of Hopf bifurcations, a curve {\tt ns} of neutral saddles, curves {\tt rhc} of rotational homoclinic bifurcations, curves {\tt snp} of saddle-node bifurcations of periodic orbits, curves {\tt chc} of contractible homoclinic bifurcations and curves {\tt sne} of saddle-node bifurcations of equilibria. The points {\tt Z}, {\tt B}, {\tt N}, {\tt X} and {\tt K} mark intersections of bifurcation curves. Reproduced with permission from
{\sc C.~Baesens and R.S.~MacKay}, {\em Resonances for weak coupling of the unfolding of a saddle-node periodic orbit with an oscillator}, Nonlinearity, 20(5):1283, 2007. \textcopyright IOP Publishing \& London Mathematical Society. All rights reserved. 
}
  \label{fig:chenciner_bifs_bae07}
\end{figure}

\Fref{fig:chenciner_bifs_bae07} reproduces from \cite{BAE07} this two-parameter bifurcation diagram of the planar vector fields on the cylinder, consisting of a number bifurcations curves and representative phase portraits in open regions and along the bifurcation curves. This figure provides a topological description of the dynamics one can to expect to find inside a Chenciner bubble, and will serve as a useful guide for what one can expect to find in the GZT~model in \sref{section:chenciner_gzt}. The phase portraits in \fref{fig:chenciner_bifs_bae07} illustrate a fundamental domain on the covering space: the left and right equilibria need to be identified in order to obtain the phase portrait on the cylinder. Hence, the saddles on the left and right sides are in fact one and the same point on the cylinder, or each others image under the $\mathbb Z_q$-symmetry. In the centre of the phase portraits are equilibria, which may be stable or unstable. The thicker ovals are contractible periodic orbits and the thicker curves that go across some phase portraits represent periodic orbits around the cylinder.

The top and bottom curves {\tt sne} in \fref{fig:chenciner_bifs_bae07} are curves of saddle-node bifurcations of equilibria, where a saddle and node are created. They are folded with respect to the vertical direction and represent the boundaries of the resonance tongue. The two curves {\tt sne} are bridged by several other bifurcation curves. The curve {\tt Hopf} connects two Takens-Bogdanov bifurcation points, denoted {\tt B}, one on each of the curves {\tt sne}. Along {\tt Hopf} the node at the centre of the phase portrait changes stability and gives birth to a contractible periodic orbit. The dashed curve {\tt ns} of neutral saddles also connects the two points {\tt B}; along it the saddle quantity, that is, the sum of the two real eigenvalues at the saddle is zero. The curve {\tt ns} does not represent a bifurcation, but it still plays a role in how the global bifurcations are organised. Emerging from the Takens-Bogdanov bifurcation points {\tt B},
are curves {\tt chc} of contractible homoclinic bifurcations, where the periodic orbit bifurcating from the curve {\tt Hopf} disappears. The curves {\tt chc} meet at a point denoted {\tt N}, where  two curves denoted {\tt rhc} intersect.

Along the curves {\tt rhc} one finds rotational homoclinic bifurcations, that is, a homoclinic connections that are non-contractible and extend around the cylindrical phase space. Each of the curves {\tt rhc} starts and ends on opposing curves {\tt sne} at points labelled {\tt Z}; these are codimension-two non-central saddle-node homoclinic bifurcations, which mark the change from a saddle-node bifurcation taking place on a periodic orbit (also referred to as SNIPER or SNIC bifurcation), or off it. The intersection point {\tt N} of the curves {\tt chc} is a codimension-two point where both non-contractible homoclinic connections occur simultaneously. The resulting homoclinic cycle can be perturbed to contractible homoclinic connections, which explains why {\tt N} is also the end point of the two curves {\tt chc}.  

The curves {\tt rhc} are intersected by the curve {\tt ns} of neutral saddles at points labelled {\tt K}, which are codimension-two bifurcations where the criticality of the non-contractible homoclinic connections changes; that is, the bifurcating non-contractible periodic orbit changes from attracting to repelling. From the points {\tt K} emanate curves  {\tt snp} of saddle-node bifurcations of periodic orbits, where a pair of attracting and repelling non-contractible periodic orbits meet and disappear. Note that one of the curves {\tt snp} also intersects the other curve {\tt rhc} at the point {\tt X}, so that three periodic orbits co-exist inside the shaded triangle formed by the points {\tt N}, {\tt K} and {\tt X}. 

The bifurcation set in \fref{fig:chenciner_bifs_bae07} represents a consistent bifurcation diagram of the dynamics of a Chenciner bubble, which is robust or structurally stable in the space of two-dimensional vector fields on a cylinder. Since \fref{fig:chenciner_bifs_bae07} takes into account all of the known bifurcations, one also speaks of a minimal bifurcation diagram of a Chenciner bubble. Since the versal unfolding of this phenomenon on the level of planar vector fields is yet unknown, it cannot be ruled out that the curves are arranged differently in the parameter plane and/or further bifurcation curves exist. This is why \fref{fig:chenciner_bifs_bae07}  represents one of several possible cases of robust bifurcation diagrams that may exist. For example, in \cite{BAE07} an extended bifurcation diagram is depicted for the case that the two Takens-Bogdanov bifurcations are oriented differently so that the curves {\tt Hopf} and {\tt ns} intersect once to form a figure-eight shape.

\subsection{Interpretation in a three-dimensional phase space}
\label{subsection:bae07_3D}

The bifurcation set and phase portraits in \fref{fig:chenciner_bifs_bae07} can be interpreted in terms of the two-dimensional Poincar\'e map $P$ and the original vector field $g$ in three-dimensions. This constitutes the inverse of the reduction to a two-dimensional approximating vector field and is a standard procedure in bifurcation analysis. The underlying idea is that the map $P$ is a perturbation of the composition of the time-one map of the normal form vector field $G$ with the action of the $\mathbb Z_q$-symmetry; for example, see \cite{GUC13,krauskopf01}. To begin with, equilibria of the phase portraits in \fref{fig:chenciner_bifs_bae07} are periodic ${p\!:\!q}$ orbits of the same stability in the three-dimensional phase space of $g$, and the curves {\tt sne} correspond to the curves of saddle-node bifurcations of periodic orbits that form the boundaries of the ${p\!:\!q}$ resonance tongue under consideration. The curve {\tt Hopf} corresponds to a curve of Neimark-Sacker or torus bifurcations of $g$, and the bifurcating contractible periodic orbit corresponds to an invariant torus around a periodic ${p\!:\!q}$ orbit. As was discussed earlier, the dynamics on the torus may be quasiperiodic or locked, and this is organised by resonance tongues emanating from this torus bifurcation curves at rational values of the frequency ratio; these secondary resonances are not described by the planar normal form $G$. 

The curves {\tt chc} and  {\tt rhc} of homoclinic bifurcations in \fref{fig:chenciner_bifs_bae07}, where the stable and unstable manifolds of the saddle equilibrium coincide, translate to a more complicated bifurcation scenario in the three-dimensional vector field $g$. Namely, the two-dimensional stable and unstable manifolds of a periodic orbit of $g$ in $\mathbb R^3$ do not coincide generically. In fact, any curve of homoclinic bifurcation of the planar normal form gives rise (under arbitrarily small generic perturbations) to a pair of curves of first and last homoclinic tangencies between these manifolds. These curves bound a wedged region in parameter space where the two manifolds intersect transversely in a homoclinic tangle. Entering this region corresponds to the break-up of the corresponding attracting or repelling torus. Here the curve {\tt rhc} corresponds to the break-up of a main torus, while the curve {\tt chc} corresponds to the break-up of a secondary torus. The respective two curves of homoclinic tangencies have exponential contact at their end points {\tt Z} and {\tt B} on the boundary curves of the ${p\!:\!q}$ resonance tongue under consideration. 

Finally, the curves {\tt snp} translate to folding tori in the three-dimensional vector field $g$. As explained above, these ``bifurcations'' of tori do not actually exist because the tori involved must break up. Thus, in the bifurcation analysis of a single Chenciner bubble one finds necessarily more folding tori with further Chenciner bubbles. More specifically, under perturbation of the time-one map of the planar normal form, the curve {\tt snp} transforms into a string of Chenciner bubbles associated with secondary (and hence weak) resonances. The result is a hierarchy of Chenciner bubbles occurring on finer and finer scales.

With this interpretation of its relevance for the full dynamics, and in spite of its lack of definiteness, \fref{fig:chenciner_bifs_bae07} provides the most comprehensive picture of the dynamics within a Chenciner bubble in a map of dimension at least two or a vector field of dimension at least three. For example, the investigation in \cite{GOV11} of a three-dimensional map, describing an adaptively controlled system, identified the (equivalent of) the bifurcation curves {\tt sne} and {\tt Hopf} connecting the boundary curves of folding resonance tongues; more specifically, the continuation software MatContM was used to find in a two-parameter plane curves of Neimark-Sacker bifurcations of periodic points that connect the respective curves of saddle-node bifurcations of periodic points bounding the respective resonance tongue. 
Very recently, the results in \cite{GOV11} were extended in \cite{neirynck18}, again with MatContM, to include the remaining bifurcation curves suggested in \cite{BAE07}.

\section{Chenciner bubbles in the GZT model}
\label{section:chenciner_gzt}

We now study the bifurcation diagram in a Chenciner bubble of the GZT~model \eref{eq:ENSO_model_GHI08}. To this end, we focus on the ${2\!:\!7}$ resonance tongue in \fref{fig:SNT_example}, because it the lowest-order (weak) resonance tongue with the most pronounced difference between the two saddle-node bifurcation curves that bound it. This helps the exposition since the respective region of the resonance tongue, which we refer to as the ${2\!:\!7}$ Chenciner bubble, is largest; moreover, working with a small period of the locked periodic solutions makes the numerical calculations more manageable. We remark that there are codimension-two Chenciner bifurcations along the curve T, where it switches between subcritical and supercritical. These points are found outside the range shown in \fref{fig:SNT_example}, for example, at $(\tau_n,c)=(0.865,3.180)$; see \cite{KEA15}.

\subsection{Bifurcation diagram in the $(\tau_n,c)$-plane}
\label{subsection:bifdiag}

\begin{figure}[t!]
  \centering
  \vspace*{7mm}
  \includegraphics{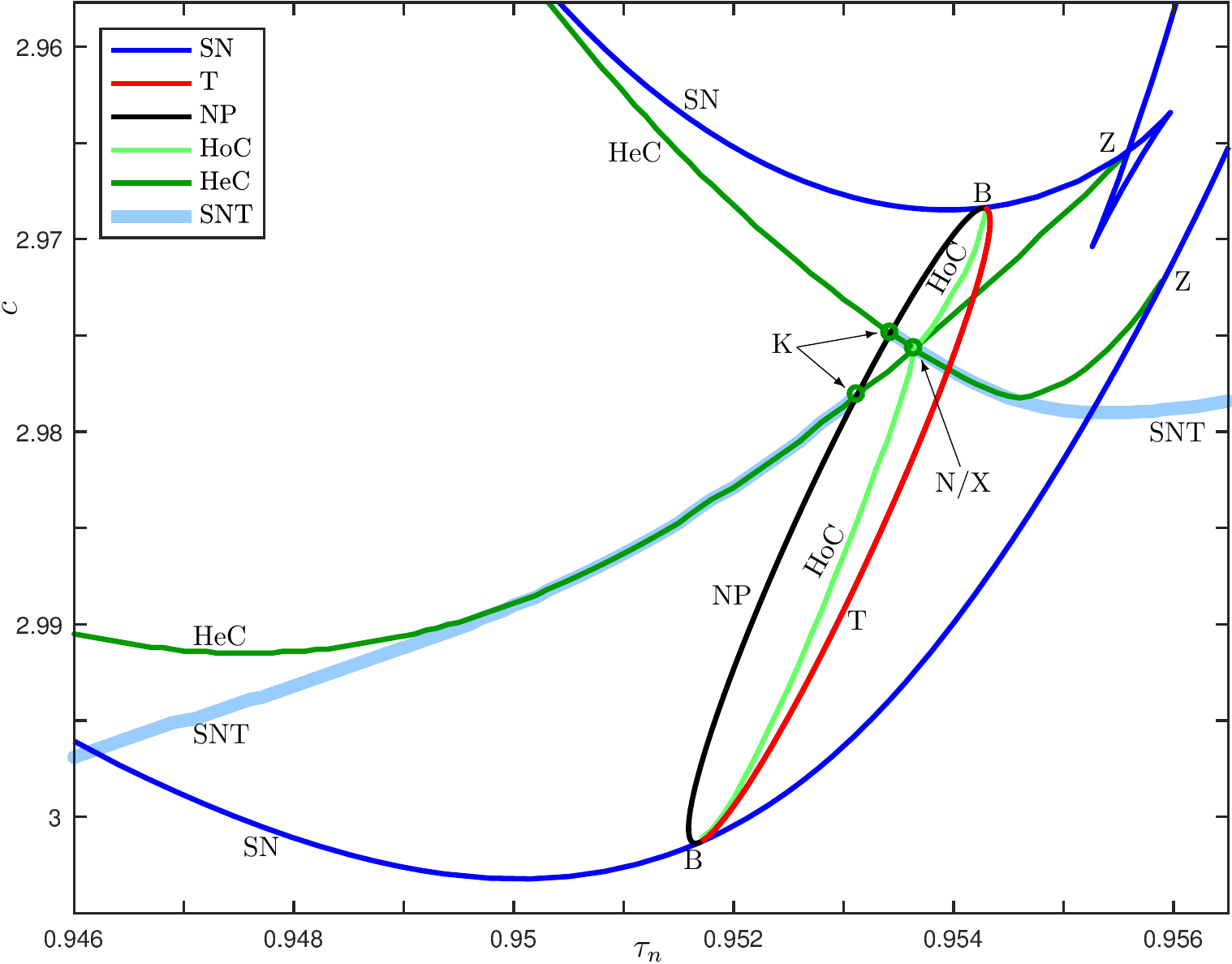}
  \caption{
Bifurcation structure of the GZT model~\eref{eq:ENSO_model_GHI08} inside the ${2\!:\!7}$ Chenciner bubble in the $(\tau_n,c)$-plane. Shown are dark blue curves SN of saddle-node bifurcations of ${2\!:\!7}$ periodic solutions, a red curve T of torus bifurcations, a black curve NP of neutral saddle periodic solutions, light green curves HoC of homoclinic transitions, dark green curves HeC of heteroclinic transitions, and light blue curves SNT of folding tori. Letters B, Z, N, X and K mark the intersection points as in \fref{fig:chenciner_bifs_bae07}. 
}
  \label{fig:chenciner_bifs}
\end{figure}

\Fref{fig:chenciner_bifs} shows the bifurcation structure inside the ${2\!:\!7}$ resonance tongue  in the $(\tau_n,c)$-plane near where it folds with respect to $c$, that is, the vertical axis as it does in \fref{fig:chenciner_bifs_bae07}. On the bounding curves SN of saddle-node bifurcations of ${2\!:\!7}$ periodic solutions in \fref{fig:chenciner_bifs} we find Bogdanov-Takens points B from which the curves T of torus bifurcations and NP of neutral saddle periodic orbits emerge; the latter is defined by $\mu_1\mu_2=1$, where $\mu_{1,2}$ are the two leading Floquet multipliers of the ${2\!:\!7}$ saddle periodic solution. The curves T and NP are calculated by continuation of periodic solutions with DDE-Biftool.  

A curve HoC of homoclinic transitions emerges from each of the two points B; the two curves HoC meet at the point N. This point N is the intersection point of two curves HeC of heteroclinic transitions, which each connect to the two opposite curves SN at points Z; note that the two points Z on the left lie outside the range shown in \fref{fig:chenciner_bifs}. From the points K where the curves HeC intersect the curve NP emerge two curves SNT of folding tori. One of the curves SNT crosses T and, hence, intersects the upper curve HoC at the point X. Because the two curves HeC and SNT are very close together, the two points N and X cannot be distinguished on the scale of \fref{fig:chenciner_bifs}. 

In this context it is important to recall that, as was discussed in \sref{subsection:bae07_3D}, the curves HoC, HeC and SNT \fref{fig:chenciner_bifs}  do not represent smooth bifurcation curves. This means that they cannot be continued with DDE-Biftool. We speak of homoclinic and heteroclinic transitions along HoC and HeC, respectively, because these represent loci of first and last tangencies, with a region of existence of a tangle in between. It turns out that the curves of first and last tangencies are so extremely close together over the range of \fref{fig:chenciner_bifs} that they cannot be distinguished for all practical purposes. The topological changes when crossing the loci HoC and HeC, on the other hand, can be detected by dedicated simulations of \eref{eq:ENSO_model_GHI08}. In this way, the locations of HoC and HeC in the $(\tau_n,c)$-plane can be determined. Specifically, we found the curves HoC by tracking the phase portrait in $\tau_n$ for a sufficient number of fixed values of $c$, and the curves HeC by tracking the phase portrait in $c$ for a sufficient number of fixed values of $\tau_n$. Similarly, the locus SNT of where tori fold is not a smooth curve, but rather the envelop of higher-order resonance tongues. It is very impractical to find these resonance tongues to determine SNT. Instead, we again employed the tracking of phase portraits in $c$ for a sufficient number of fixed values of $\tau_n$ to determine the two loci SNT; as one may expect for resonances of very high order, the loci SNT in \fref{fig:chenciner_bifs} are effectively one-dimensional. We remark that determining the loci HoC, HeC and SNT in the $(\tau_n,c)$-plane is a considerable effort. 

Overall, we find very good agreement between the computed bifurcation diagram \fref{fig:chenciner_bifs} of the GZT model~\eref{eq:ENSO_model_GHI08} and the suggested bifurcation diagram \fref{fig:chenciner_bifs_bae07} from \cite{BAE07}. Indeed, the different curves of bifurcations and transitions agree, taking into account the interpretation of the planar normal form for the dynamics of periodic orbits and tori of the DDE. Examples of phase portraits and the respective transitions will be shown and discussed in \sref{subsection:T_HoC}--\ref{subsection:HeC_SNT}, where we demonstrate how the dynamics of the GZT model~\eref{eq:ENSO_model_GHI08} changes when the bifurcation diagram in \fref{fig:chenciner_bifs} is crossed.

\subsection{Criticality of the curve T}
\label{subsection:Tcrit}

We now discuss briefly a feature of the bifurcation diagram \fref{fig:chenciner_bifs} that differs from the bifurcation diagram in \fref{fig:chenciner_bifs_bae07}, namely the criticality of the curve T. In \cite{BAE07} it is assumed that the Hopf bifurcation is supercritical throughout, so that a family of stable limit cycles exists between the curve {\tt Hopf} and the curve {\tt chc}. However, simulations of the GZT~model~\eref{eq:ENSO_model_GHI08} reveal no evidence of a stable torus being born along curve T, suggesting that the torus bifurcation is subcritical throughout. To determine its criticality we now compute resonance tongues of locked periodic solutions that emerge from curve T. 

\begin{figure}[t!]
  \centering
  \vspace*{5mm}
  \includegraphics{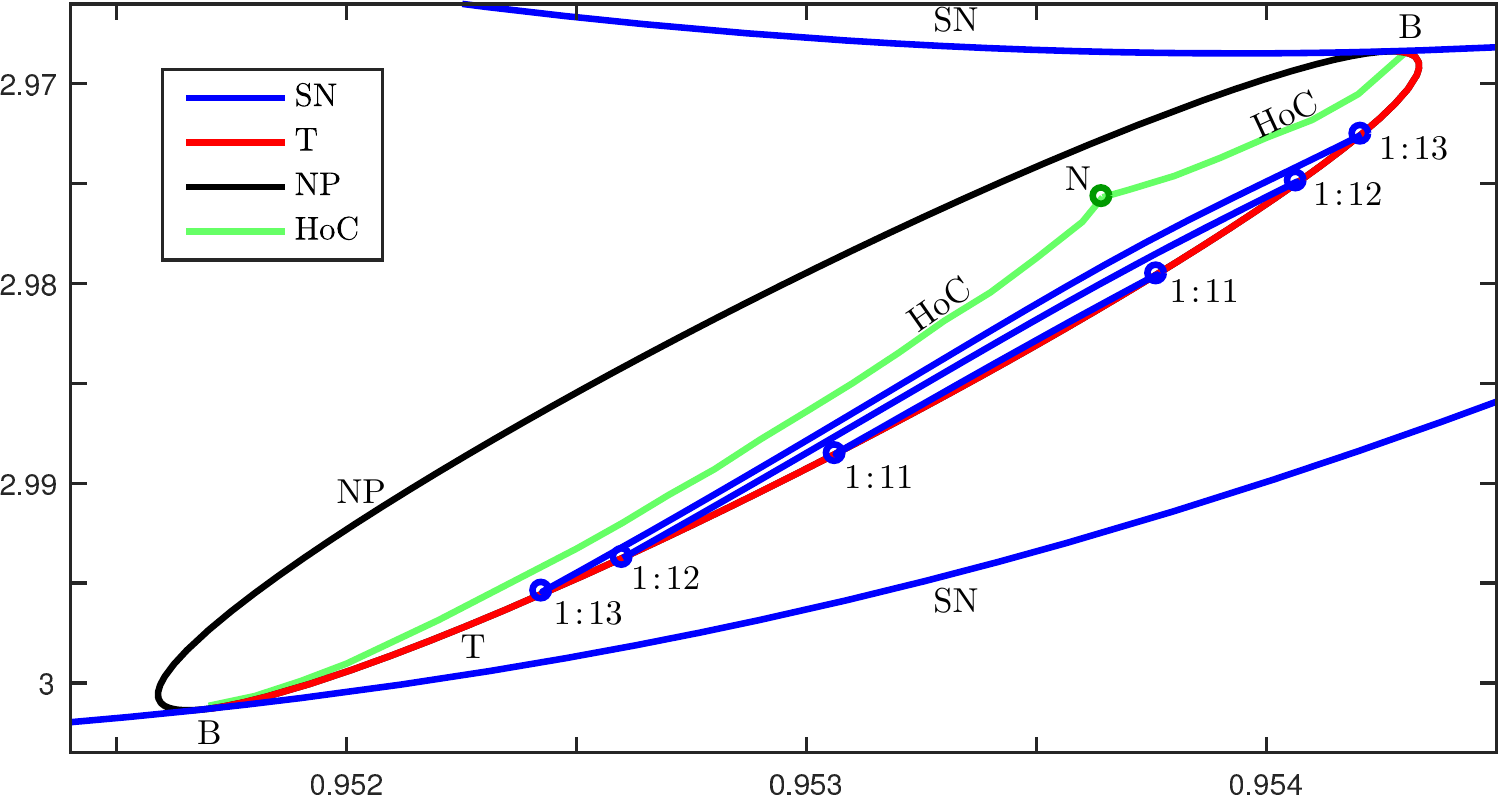}
  \caption{
Enlarged view of the bifurcation diagram from \fref{fig:chenciner_bifs}. Shown are only the curves SN, T, NP and HoC, and also blue resonance tongues emerging at resonance points from T; the resonances are labelled ${1\!:\!q}$, which means that the locked periodic solution makes $q$ turns around the ${2\!:\!7}$ periodic solution. 
}
  \label{fig:T_critical}
\end{figure}

\Fref{fig:T_critical} is an enlargement of the bifurcation diagram \fref{fig:chenciner_bifs}, where only curves directly related to the criticality of the curve T are displayed. Again, the two curves SN at the very top and bottom of \fref{fig:T_critical} are the boundaries of the original ${2\!:\!7}$ resonance tongue. The additional blue curves represent secondary resonance tongues that are rooted on the curve T, at points of ${1\!:\!11}$, ${1\!:\!12}$ and ${1\!:\!13}$ resonance with respect to the ${2\!:\!7}$ periodic solution. Hence, the locked periodic solutions in the resonance tongues have periods $77$, $84$ and $91$ with respect to the one year forcing cycle, respectively. Because they are extremely narrow, as we determined numerically, we compute these resonance tongues via the continuation of a single curves of the respective locked periodic solution. The secondary resonance tongues are in the region where the ${2\!:\!7}$ periodic solution is stable. Moreover, the secondary locked periodic solutions have one unstable Floquet multiplier, which means that they exist on a torus of saddle-type. This confirms that T is indeed a curve of subcritical torus bifurcations. As we have found and will discuss next, the tori of saddle-type that are created at the curve T terminate at curve HoC. This is in contrast to the suggested bifurcation diagram in \fref{fig:chenciner_bifs_bae07}, where the curves {\tt Hopf} and {\tt chc} of contractible homoclinic bifurcations both concern stable periodic orbits.

\subsection{Transition through the curves T and HoC}
\label{subsection:T_HoC}

\begin{figure}[p]
  \centering
  \includegraphics[width=\textwidth]{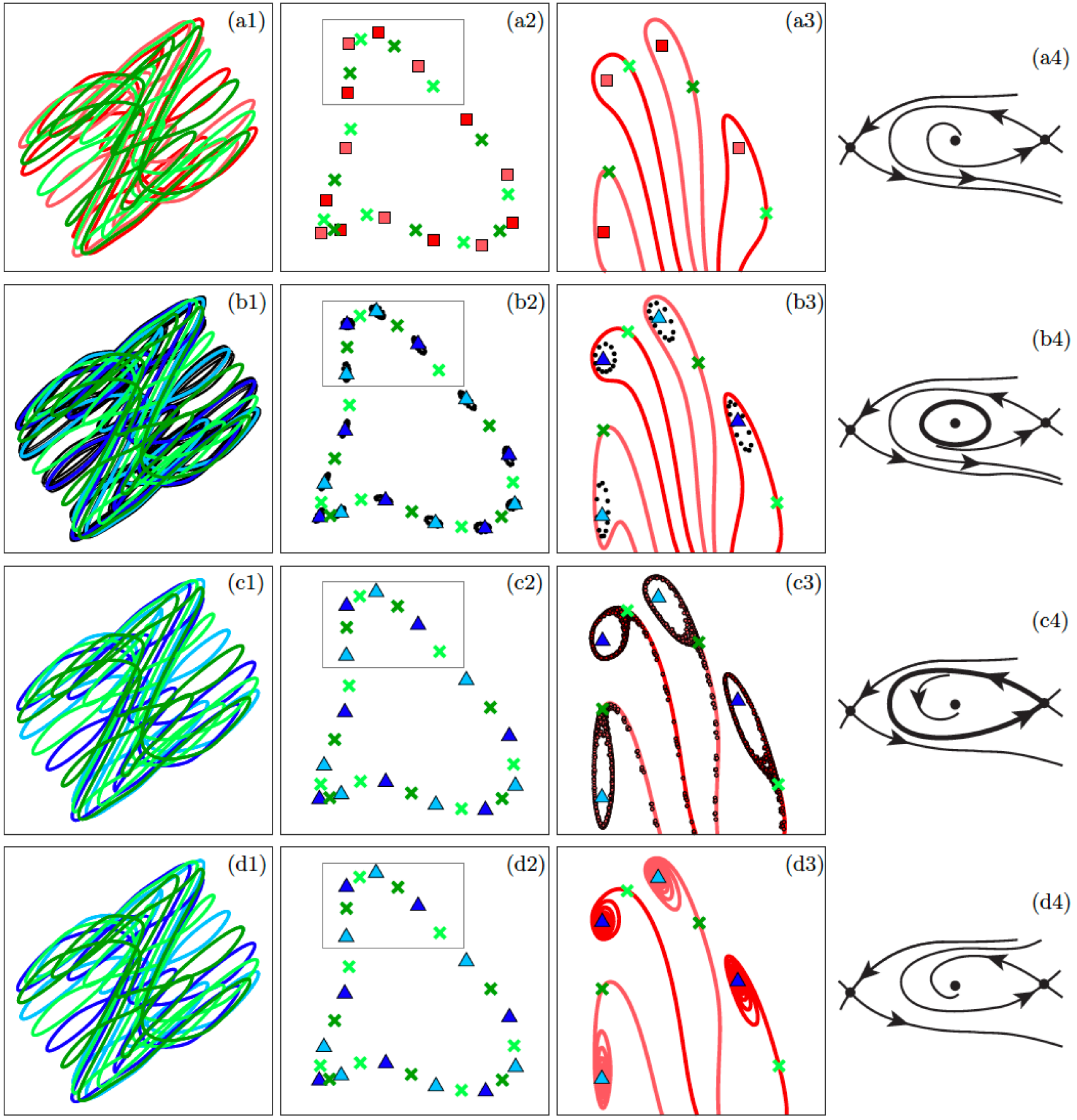}
  \caption{Phase portraits in the transition through the curves T and HoC, for $c=2.9850$ and $\tau_n=0.9540$ (a), $\tau_n=0.9533$ (b), $\tau_n=0.9531$ (c) and $\tau_n=0.9530$ (d).
shown are  projections onto the $(h(t),h(t-\tau_n))$-plane of periodic solutions (blue, green and red curves) and a saddle torus (black) (first column); the respective stroboscopic trace in the $(h(t),h(t-\tau_n))$-plane (second column); an enlargement of the grey box that also shows the unstable manifolds (red curves and dots in panel~(c3)) of the 1-saddles (third column); and a sketch of the corresponding phase portrait of the planar normal form (fourth column). Attracting periodic solutions are blue, 1-saddles are green and 2-saddles are red; the two different symmetry-related periodic solutions are distinguished by different shades of the respective colour. 
}
  \label{fig:T_HoC_examples}
\end{figure}

\Fref{fig:T_HoC_examples} shows how the phase portrait of the GZT model~\eref{eq:ENSO_model_GHI08} changes, for  $c=2.9850$ and decreasing $\tau_n$, during the transition through the curves T and HoC. This figure also serves to illustrate how the curve HoC was identified in the DDE. More specifically, we compute periodic solutions and tori and present them in different ways: in projection onto the $(h(t),h(t-\tau_n))$-plane and as stroboscopic trace in the $(h(t),h(t-\tau_n))$-plane, where triangles, crosses and squares represent stable, 1-saddle and 2-saddle periodic solutions, respectively. An enlargement of the stroboscopic trace also shows the one-dimensional trace of the unstable manifolds of the 1-saddle periodic solutions; as in \cite{CAL16} these curves were calculated via integration of initial conditions very close to the 1-saddle periodic solutions along their unstable eigendirections --- technically, within a region referred to as a fundamental domain \cite{KRA03}. Throughout, light and dark shades of the same colours distinguish between the two symmetry-related periodic solutions (whose existence was  discussed in \sref{section:intro}); see \cite{KEA15} for more details.  To facilitate the comparison, we provide the corresponding sketch of the normal form on the cylinder in the style of \fref{fig:chenciner_bifs_bae07}. 

For parameters $\tau_n$ to the right of curve T in \fref{fig:chenciner_bifs}, as in row~(a) of \fref{fig:T_HoC_examples}, there exists two symmetry-related 1-saddle and two symmetrically-related 2-saddle ${2\!:\!7}$ periodic solutions; see panel (a1). In the stroboscopic trace in the $(h(t),h(t-\tau_n))$-plane of panels (a2) and (a3) they appear as two sets each of period-seven points. The traces of the unstable manifolds of the 1-saddles surround the 2-saddles and converge to a small-amplitude stable torus (not shown in \fref{fig:T_HoC_examples} for clarity and corresponding to the green/yellow region of the maximum maps in \fref{fig:SNT_example}). When the curve T in \fref{fig:chenciner_bifs} is crossed for decreasing $\tau_n$, the 2-saddle periodic solutions become stable and a torus is created. This is illustrated in row~(b) of \fref{fig:T_HoC_examples}. Notice how the torus, which we find as a high-period locked periodic solution, surrounds the respective stable periodic solution; see, in particular, panel (b3) and notice that the unstable manifolds of the 1-saddles effectively remain unchanged. 

For $\tau_n=0.9531$, shown in row~(c) of \fref{fig:T_HoC_examples}, there is evidence that the unstable manifolds form a homoclinic tangle. 
More specifically, the unstable manifolds in one direction still converge to the smaller torus. However, the unstable manifolds in the other direction have become more complicated. After encircling the attracting periodic points, some of the integrated points from the fundamental domain converge to the smaller torus, some converge towards the attracting periodic points and some accumulate at the 1-saddle points. In fact, this distribution of integrated points corresponds to an unstable manifold that is subject to infinite folding and stretching, as is typical for a homoclinic tangle. Because of the drastic folding and stretching (of the trace) of the unstable manifold, we do not represent the unstable manifold as curves; this would require considerably more sophisticated calculations \cite{KRA03}.
For the practical purpose of detecting the locus HoC in the bifurcation diagram in \fref{fig:chenciner_bifs}, the representation of the homoclinic tangle in panel~(c3) of \fref{fig:T_HoC_examples} provides a clear identification of the homoclinic transition.
Finally, row~(d) shows the situation to the left of HoC, when the unstable manifolds of the 1-saddles have undergone the homoclinic transition HoC and now spiral onto the attracting periodic solution. 
Note that we do not find evidence specifically for the first and last homoclinic tangencies and conclude that they occur within a very small parameter range on the scale of the ${2\!:\!7}$ Chenciner bubble.

\Fref{fig:T_HoC_examples} clearly shows that crossing the curves T and HoC in \fref{fig:chenciner_bifs} produces qualitative changes of the dynamics that correspond to crossing the curves {\tt Hopf} and {\tt chc} in \fref{fig:chenciner_bifs_bae07}, when taking into account that the curve T is now subcritical; compare panels (a4)--(c4) of \fref{fig:T_HoC_examples} with the respective phase portraits of \fref{fig:chenciner_bifs_bae07}.

\subsection{Transition through the curves HeC and SNT}
\label{subsection:HeC_SNT}

\begin{figure}[p]
  \centering
  \includegraphics{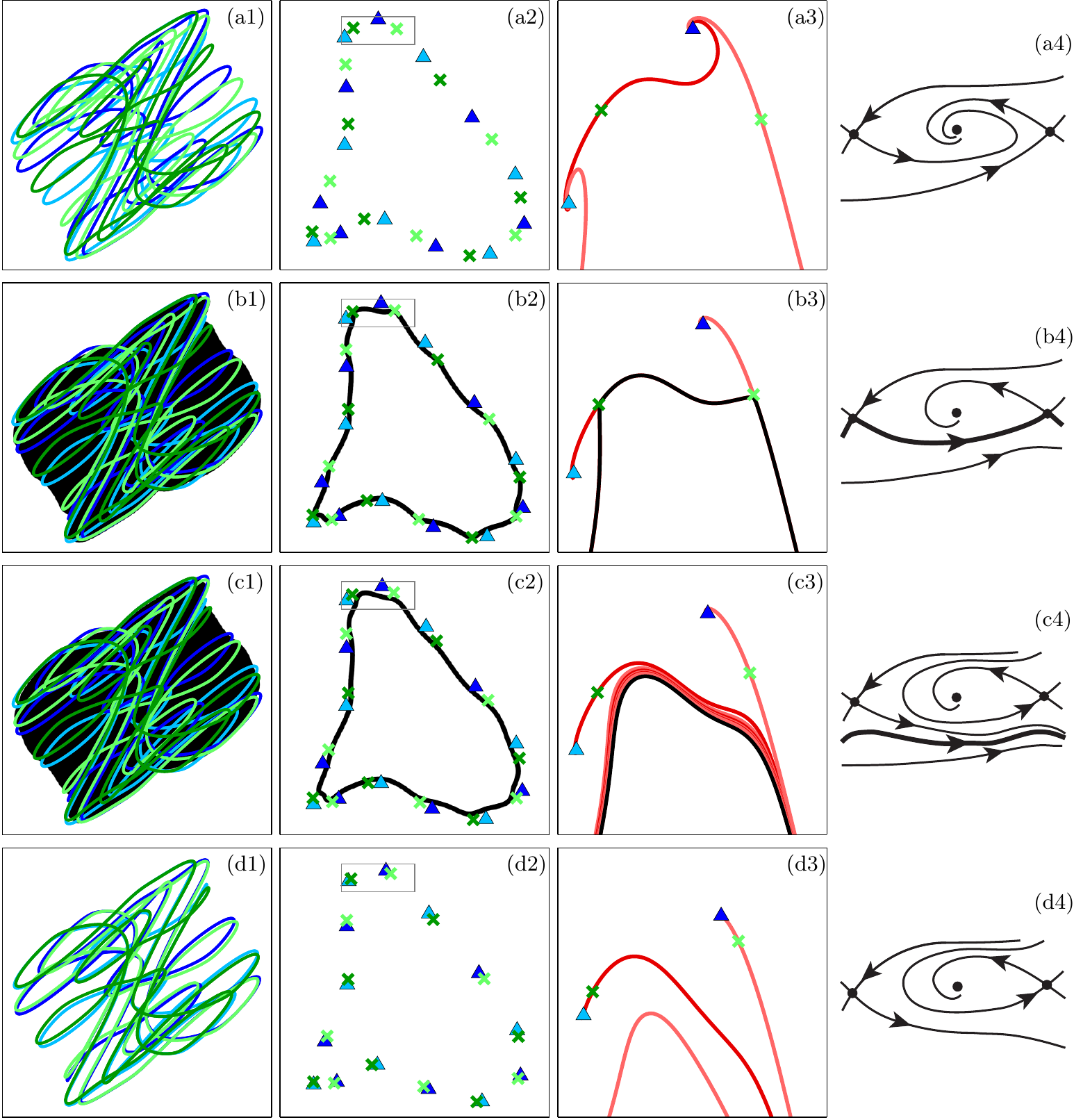}
  \caption{Phase portraits in the transition through the curves HeC and SNT, for  $\tau_n=0.9475$ and $c=2.9850$ (a), $c=2.9915$ (b), $c=2.9940$ (c) and $c=2.9990$ (d).
Shown are projections onto the $(h(t),h(t-\tau_n))$-plane of periodic solutions (blue and green curves) and a saddle torus (black) (first column); the respective stroboscopic trace in the $(h(t),h(t-\tau_n))$-plane (second column); an enlargement of the grey box that also shows the unstable manifolds (red curves) of the 1-saddles (third column); and a sketch of the corresponding phase portrait of the planar normal form (fourth column). Attracting periodic solutions are blue and 1-saddles are green; the two different symmetry-related periodic solutions are distinguished by different shades of the respective colour. 
}
  \label{fig:HeC_SNT_examples}
\end{figure}

\Fref{fig:HeC_SNT_examples} shows similarly how the phase portrait of the GZT model~\eref{eq:ENSO_model_GHI08} changes during the transition through the curves HeC and SNT for $\tau_n=0.9475$ and increasing $c$. Row~(a) is for parameter $c$ above both curves HeC and SNT in \fref{fig:chenciner_bifs}, where there exists a pair of symmetry-related stable periodic solutions and a pair of symmetry-related 1-saddle periodic solutions; see \fref{fig:HeC_SNT_examples}(a1) and (a2). As panel (a3) shows, the two branches of the unstable manifolds of the 1-saddles spiral into neighbouring stable points in the trace. In particular, the spiralling means that the attracting torus is not normal hyperbolic.  The moment of the heteroclinic transition HeC is illustrated in row~(b) of \fref{fig:HeC_SNT_examples}: one branch of the unstable manifold of each 1-saddle connects to the respective other symmetry-related 1-saddle. The black curve that illustrates these connections in panel (b3) is actually a stable torus; see also panels (b1) and (b2). Since this invariant curve passes exceptionally close to the periodic points, this situation can be considered the moment of the HeC transition for the practical detection of the curve HeC in the bifurcation diagram. In contrast to the HoC transition in \sref{subsection:T_HoC}, we found no evidence of first and last tangencies or tangles. In between the curves HeC and SNT, as in row~(c) of \fref{fig:HeC_SNT_examples} the torus is smooth and further away from the 1-saddles; one branch of their unstable manifolds accumulates on the torus, while the other branch goes to the nearby attracting periodic solution. At the curve SNT the stable torus disappears by colliding with a saddle torus. Even though the locus SNT is associated with secondary and higher-order Chenciner bubbles, numerically, this appears to happen instantly and can be detected as such to find SNT as a curve in \fref{fig:chenciner_bifs}. The phase portrait is then as shown in row~(d) of \fref{fig:HeC_SNT_examples}. Again, the sketches of the corresponding phase portraits of the planar normal form clearly show that the transition through HeC and SNT in \fref{fig:chenciner_bifs} is described well by that through the curves {\tt rhc} and {\tt snp} in \fref{fig:chenciner_bifs_bae07}.

\newpage

\section{Discussion and interpretation for tipping}
\label{section:disc}

This paper took a detailed look at folding tori in a DDE model, which represents the interaction between negative delayed feedback and periodic forcing. In particular, this phenomenon involves Chenciner bubbles, where the tori break up and the dynamics becomes complicated. We analysed the dynamics inside the ${2\!:\!7}$ Chenciner bubble of the GZT~model by computing curves of bifurcations of periodic solutions and tori by means of continuation and numerical integration; for some curves this required the computation of unstable manifolds of saddle periodic solutions. The bifurcation set we found inside the ${2\!:\!7}$ Chenciner bubble compares very well with the suggested bifurcation diagram from \cite{BAE07} of a planar normal form. 

Folding tori are of particular interest in the context of the phenomenon of (climate or other) tipping. One speaks of a tipping event when a slight parameter variation of the system creates a qualitative or drastic response \cite{LEN08}. Such sudden changes in observed behaviour can be associated with certain bifurcations \cite{THO11} . For example, past data and various mathematical models, such as those in \cite{LEN09,MAN88,MAR04,RAH95,RAH05}, provide evidence of bistability with a fold bifurcation in the North Atlantic thermohaline circulation (THC) --- the circulation mechanism that drives the Gulf Stream. In past studies only simple bifurcations of equilibria or periodic solutions are considered. Indeed, the standard case of tipping is that of an attracting equilibrium or periodic orbit approaching a fold bifurcation, in which it then disappears. We propose that the passage of a stable torus through the locus of folding tori is also a type of tipping. In light of the fact that at least two frequencies need to be involved, we refer to this kind of tipping as multi-frequency tipping. 

The additional complexity of the dynamics we encountered in folding tori may have important implications for tipping events in systems that are driven by and/or feature multiple frequencies. 
This is certainly the case for climate systems more generally, which are driven by effects including orbital forcing at various time scales, as well as different types of feedback mechanisms; for example, see \cite{stocker01} and references therein. 
A critical and intriguing question is whether multi-frequency tipping can be identified in real-world data or 
models of at least intermediate complexity. In fact, there are already some possible candidates for multi-frequency tipping events. For example, the above mentioned North Atlantic THC is known to exhibit fold bifurcations. In \cite{HEL04} a time series analysis is used to predict a tipping event induced by an unspecified bifurcation in a THC model of intermediate complexity. Due to the many degrees of freedom in the model, the detected bifurcation may, in fact, be a quasiperiodic bifurcation involving tori. As a further example, time series analysis is used in \cite{DAK08} to identify, in data taken from tropical Pacific sediment cores, an ancient greenhouse to icehouse tipping event --- described as ``fold-like'' in \cite{THO11}. Given the many degrees of freedom of real-world climate systems, this tipping event may be most accurately described by folding tori as well.

A particular aspect of multi-frequency tipping is that Chenciner bubbles are crossed which, depending on the frequencies involved, may be sufficiently large to be of relevance. This raises a number of interesting questions: What influence does the additional bifurcation structure inside the Chenciner bubble have on the observed dynamics when tipping is approached? Can this knowledge provide useful precursors to warn of an imminent multi-frequency tipping event? 
While answering these question will require considerable follow-on work, we now provide a brief example of multi-frequency tipping that highlights the role of the relative time scales involved. 

\begin{figure}[t]
  \centering
  \includegraphics{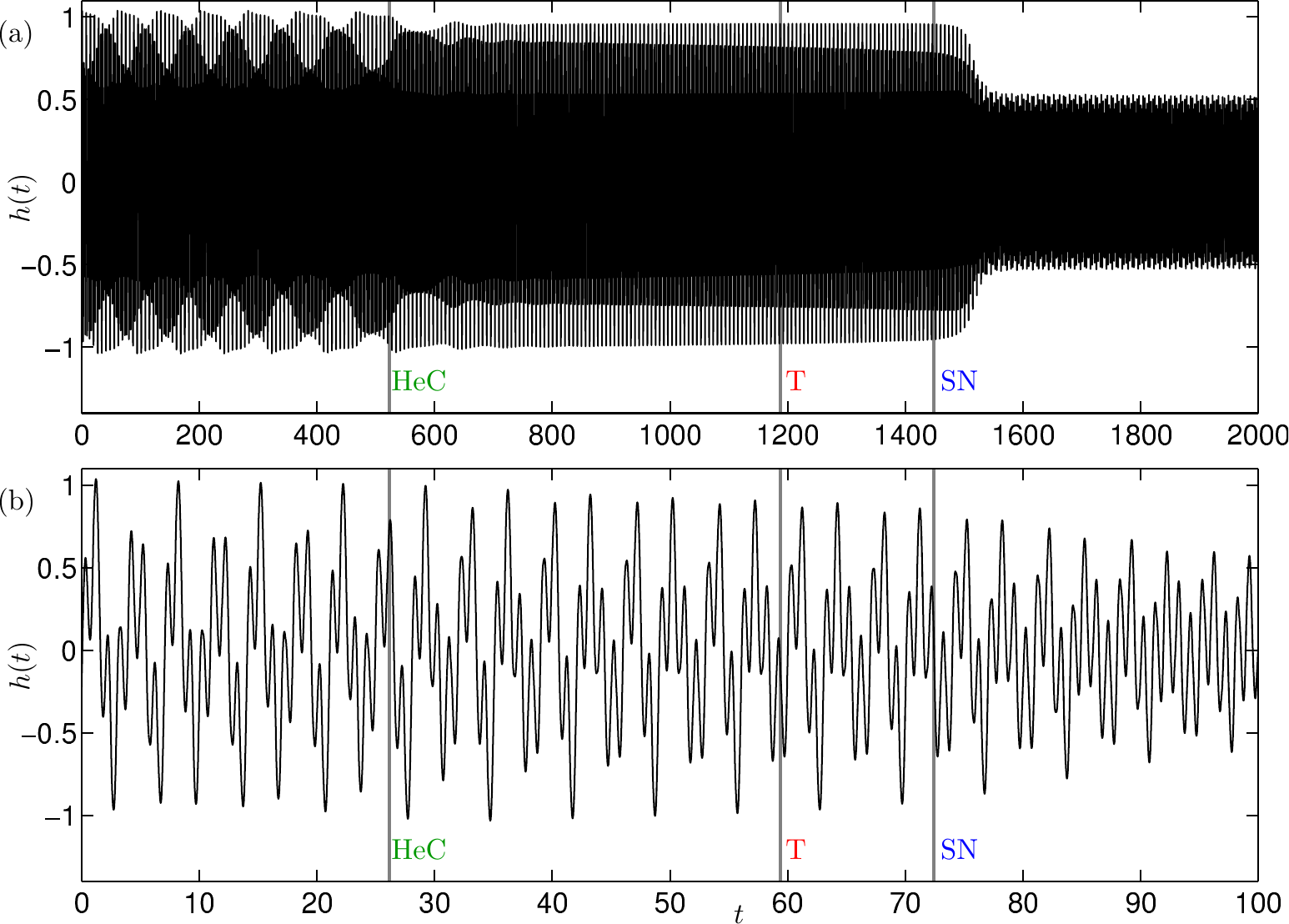}
  \caption{
Times series of the GZT~model for fixed $\tau_n=0.953$  but with $c$ drifting linearly, from 2.96 to 3.01, across the ${2\!:\!7}$ Chenciner bubble, over 2000 years in panel~(a) and over 100 years in panel ~(b). Grey vertical lines indicate where $c$ crosses the curves HeC, T and SN in \fref{fig:chenciner_bifs}.
}
  \label{fig:drifting_c}
\end{figure}

\Fref{fig:drifting_c} shows two time series of the GZT~model~\eref{eq:ENSO_model_GHI08} for $\tau_n=0.953$ fixed and with $c$ drifting linearly from 2.96 to 3.01, that is, across the ${2\!:\!7}$ Chenciner bubble shown in \fref{fig:chenciner_bifs}. The time series in panels~(a) and~(b) of \fref{fig:drifting_c} drift across the Chenciner bubble in 2000 and in 100 years, respectively. While doing so, the curves HeC, T and SN are crossed, and this is indicated by the labelled vertical lines in the time series. This corresponds to the transition from multi-frequency dynamics on a large torus to dynamics on a considerably smaller torus.

In \fref{fig:drifting_c}(a), for a very slow drift in $c$, one can clearly associate the observed behaviour in the time series with the quasi-static passage through the different bifurcations. Before HeC is reached, the dynamics is characterised by an amplitude modulation of about 40 years. After passing through HeC, at $t \approx 530$, the time series appear to settle on the stable ${2\!:\!7}$ periodic solution, which results in a much shorter period of the amplitude modulation. When T is crossed, at $t \approx 1200$, the ${2\!:\!7}$ periodic solution loses stability in a torus bifurcation, but there is no obvious change in the time series because the instability is only very weak. Finally, the boundary SN of the folding ${2\!:\!7}$ resonance tongue is crossed at $t \approx 1450$. As a result, the ${2\!:\!7}$ periodic solution ceases to exist and the time series settles onto quasiperiodic oscillations of a much smaller amplitude. One could argue that, in this scenario, the heteroclinic transition HeC acts as a precursor to the more dramatic transition at SN. In \fref{fig:drifting_c}(b) the drift occurs over only 100 years. There is a notable change of the frequency of the observed oscillation at HeC, which may again be interpreted as a precursor for the transition to much smaller amplitude oscillations when SN is crossed. However, the changes in the behaviour of the time series appear to be more subtle for this faster drift. 

A wider understanding of multi-frequency tipping may be of interest in the contexts of both long-term (palaeoclimatology) and short-term (anthropogenic or human-induced) climate change, and likely also in other fields, such as ecology. In ongoing and future work, it will be very interesting to study specific candidate systems with folding tori to try and identify qualitative changes in their time series data that can be associated with Chenciner bubbles. 

Asides from the effect of different time scales, or rather the size of the relevant Chenciner bubble in parameter space, there are other factors that could impede the identification of torus break-up in time series data. 
Parameter uncertainty could be an issue. In the brief case study presented here, we have analysed only a single Chenciner bubble that exists for a small range of parameters. However, the ${2\!:\!7}$ Chenciner bubble presented here is only one of an infinite number of Chenciner bubbles, of varying sizes, that exist along the locus of folding tori. For example, in \fref{fig:SNT_example} one would expect a string of Chenciner bubbles along the sharp interface between large and small maxima.
Moreover, there are of course many different parameter paths through a Chenciner bubble that a quasi-static solution could take, and this will affect the specific dynamics that are encountered. 

A further issue in this context is that time series data taken from observations or complicated models will be subject to noise, for example, reflecting high-frequency effects omitted in conceptual models such as the GZT~model. Understanding the extent to which noise can conceal the effects of the bifurcation structure of a Chenciner bubble in data is ongoing work.
Finally, more complicated models, compared to the GZT~model considered here, will generally dependent not only on time but also on space. Instead of analysing time series data for an entire mesh of points in space, the data is generally converted into an observable (for example, mean sea-surface temperature across the Pacific Ocean) or reduced to a small number of critical components by, for example, a principal component analysis in order to identify strong trends in the data. It is an interesting and practical challenge to determine to what extent such averaging or data reductions preserve signatures of torus break-up.

\subsection{Acknowledgements}
We thank Claire Postlethwaite and Anup Purewal for helpful discussions and Jan Sieber for his assistance with DDE-Biftool, especially regarding the code for the continuation of neutral saddles. We acknowledge the New Zealand eScience Infrastructure for use of their high-performance computing cluster.

\bibliographystyle{siam}
\bibliography{kk_bubbles_refs}


\end{document}